\documentclass[11pt]{article}
\usepackage{latexsym,amssymb,amsmath,a4wide,theorem}
\usepackage[isolatin]{inputenc}
\usepackage[T1]{fontenc}
\usepackage{times}
\usepackage{cite}
\usepackage{amsfonts}
\usepackage{hyperref}

\def\sqr#1#2{{\vbox{\hrule height.#2pt
     \hbox{\vrule width.#2pt height#1pt \kern#1pt
           \vrule width.#2pt}
     \hrule height.#2pt}}}

\newtheorem{theorem}{Theorem}[section]
\newtheorem{lemma}[theorem]{Lemma}
\newtheorem{remark}[theorem]{Remark}

\newtheorem{corollary}[theorem]{Corollary}

\title{Multiplicity of solutions for a class of critical Schr\"odinger-Poisson system with two parameters}

\author{
{Yongpeng Chen$^{1}$},
{Zhipeng Yang$^{2}$}\thanks{Corresponding author: zhipeng.yang@mathematik.uni-goettingen.de}\\
\small School of Science, Guangxi University of Science and Technology, Liuzhou 545006. P.R.China.$^{1}$\\
\small Mathematical Institute, Georg-August-University of G\"ottingen, G\"ottingen 37073, Germany.$^{2}$\\
}

\date{}

\begin{document}
\maketitle
\begin{abstract}
We study a class of critical Schr\"odinger-Poisson system of the form
\begin{equation*}
\begin{cases}
-\Delta u+\lambda V(x)u+\phi u=\mu |u|^{p-2}u+|u|^{4}u& \quad x\in \mathbb{R}^3,\\
-\Delta \phi=u^2&\quad x\in \mathbb{R}^3,\\
\end{cases}
\end{equation*}
where $\lambda, \mu>0$ are two parameters, $p\in(4,6)$ and $V$ satisfies some potential well conditions. By using the variational arguments, we prove the existence of positive ground state solutions for $\lambda$ large enough and $\mu>0$, and their  asymptotical behavior  as $\lambda\to\infty$. Moreover, by using  Ljusternik-Schnirelmann theory, we obtain the existence of multiple positive solutions if $\lambda$ is large and $\mu$ is small.
\end{abstract}


{\bf AMS} Subject Classification: 35A15, 35B40, 35J20
\numberwithin{equation}{section}
\section{Introduction and Main Results}

In this paper, we consider the following critical Schr\"odinger-Poisson system
\begin{equation}\label{maineq}
\begin{cases}
-\Delta u+\lambda V(x)u+\phi u=\mu |u|^{p-2}u+|u|^{4}u& \quad x\in \mathbb{R}^3,\\
-\Delta \phi=u^2&\quad x\in \mathbb{R}^3,\\
\end{cases}
\end{equation}
where $\lambda, \mu>0$ are two parameters, $p\in(4,6)$ and the potential $V(x)$ is a locally H\"oder continuous function satisfies the following assumptions:
\begin{itemize}
\item[($V_1$)] $V(x)\geq0$ on $\mathbb{R}^3$.
\item[($V_2$)] There exists $M_0>0$ such that the set $A:=\{x\in\mathbb{R}^3:a(x)\leq M_0\}$ is nonempty and $m(A)<\infty$, where $m(A)$ denotes the Lebesge measure of $A$ on $\mathbb{R}^3.$
\item[($V_3$)] $\Omega:=int\{V^{-1}(0)\}$ is a nonempty smooth bounded  domain and $\bar{\Omega}=V^{-1}(0)$. Without loss of generality, we assume $0\in\Omega$.
\end{itemize}
\par
We are interested in the existence of positive ground state solutions for $\lambda$ big enough and $\mu>0$, and their  asymptotical behavior  as $\lambda\to\infty.$ Especially, the existence of multiple positive solutions is investgated when $\lambda$ is large and $\mu$ is small.
The hypothesis about the potential $V(x)$ were first introduced by Bartsch and Wang \cite{Bartsch-Wang95CPDE} in the study of a nonlinear Schr\"{o}dinger equation and has been widely used by many authors, for instance, see \cite{Bartsch-Wang01CCM,Ding-Szulkin07CVPDE,Jiang-Zhou11JDE,Sun-Wu14JDE,Ye-Tang15CVPDE,Zhao-Liu-Zhao2013JDE} and the reference therein. The conditions $(V_1)-(V_2)$ imply that $\lambda V$ represents a potential well whose depth is controlled by $\lambda$. Therefore, $\lambda V$ is referred as the steep potential well if $\lambda$ is sufficiently large and one expects to find solutions which localize near its bottom $\Omega$.
\par
In recent years, a great deal of work has been devoted to the study of standing waves for the Schr\"{o}dinger equation:
\begin{equation}\label{eq1.2}
 i\frac{\partial \Psi}{\partial t}=-\Delta\Psi+\widetilde{V}(x)\Psi-f(|\Psi|)\ \ \ \ \text{in}\ \mathbb{R}^3 \times \mathbb{R},
 \end{equation}
where $i$ is the imaginary unit, $\widetilde{V}(x)=V(x)+E$ is the potential function with the constant $E$ and $f(exp(i\theta)\xi)=exp(i\theta)f(\xi)$ for $\theta,\xi\in\mathbb{R}$ is a nonlinear function.
\par
An interesting Schr\"{o}dinger equation class is when it describes quantum (nonrelativistic) particles interacting with the electromagnetic field generated by the motion. That is a nonlinear Schr\"{o}dinger-Poisson system (also called  Schr\"{o}dinger-Maxwell system):
\begin{equation}\label{eq1.3}
\begin{cases}
i\frac{\partial \Psi}{\partial t}=-\Delta\Psi+\widetilde{V}(x)\Psi+\phi\Psi-f(|\Psi|)& \text{in}\  \mathbb{R}^3\times\mathbb{R},\\
-\Delta\phi=|\Psi|^2& \text{in}\ \mathbb{R}^3.
\end{cases}
\end{equation}
A solution of the form $(e^{-iEt}u(x), \phi(x))$ is called a standing wave and $(e^{-iEt}u(x), \phi(x))$ is a solution of \eqref{eq1.3} if and only if $(u(x),\phi(x))$ satisfies
\begin{equation}\label{eq1.5}
\begin{cases}
-\Delta u+V(x)u+\phi u=f(u)& \text{in}\  \mathbb{R}^3,\\
-\Delta\phi=u^2& \text{in}\ \mathbb{R}^3,
\end{cases}
\end{equation}
which was proposed by Benci and Fortunato \cite{Benci-Fortunato1998TMNA} in 1998 on a bounded domain,  and is related to the Hartree equation (\cite{Lions1987}). Recently, in order to better simulate the interaction effect among many particles in quantum mechanics, systems of a nonlinear version of the Schr\"{o}dinger equation coupled with a Poisson equation have begun to receive much attention. We refer the interested readers to see \cite{Ambrosetti2008MJM,Azzollini2010JDE,Azzollini-d'Avenia-Pomponio2010AIHPAN,Cerami-Vaira2010JDE,Jiang-Zhou11JDE,Zhao-Liu-Zhao2013JDE} and the references therein.
\par
Recently, Jiang and Zhou \cite{Jiang-Zhou11JDE} first applied the steep potential well conditions $(V_1)-(V_3)$ to Schr\"{o}dinger-Poisson system, and proved the existence of solution. Moreover, they also studied the asymptotic behavior of solution by combining domains approximation with priori estimates. Later, Zhao et al. \cite{Zhao-Liu-Zhao2013JDE} considered a case allowing the potential $V$ changes sign. That is, for the following system
\begin{equation*}
\begin{cases}
-\Delta u+\lambda V(x)u+K(x)\phi u=|u|^{p-2}u& \text{in}\ \mathbb{R}^3,\\
-\Delta \phi=K(x)u^2& \text{in}\ \mathbb{R}^3,
\end{cases}
\end{equation*}
the authors assume that $V\in C(\mathbb{R}^3,\mathbb{R})$ is bounded from below, and satisfies $(V_2)-(V_3)$.
Using variational setting of \cite{Ding-Wei07JFA}, they obtained the existence and asymptotic behavior of nontrivial solutions for $p\in(3,6)$.
\par
In \cite{MR2107669}, Clapp and Ding have studied the nonlinear Schr\"{o}dinger equation
\begin{equation*}
-\Delta u+\lambda V(x) u=\mu u+u^{2^{*}-1}, \quad \text { in } \mathbb{R}^{N}
\end{equation*}
for $N \geq 4, \lambda, \mu>0$ and $V$ verifying $\left(V_{1}\right)-\left(V_{3}\right) .$ By using variational methods, the authors established existence and multiplicity of positive solutions which localize near the potential well for  $\lambda$ large and $\mu$ small. Later, Alves and Barros \cite{MR3850308} generalized these results to
\begin{equation}\label{eq1.6}
-\Delta u+\lambda V(x) u=\mu u^{p-1}+u^{2^{*}-1}, \quad x \in \mathbb{R}^{N}.
\end{equation}
By combining variational methods with the Ljusternik-Schnirelmann category, the authors were able to show that there exist $\lambda^{*}, \mu^{*}>0$ such that for $\lambda \geq \lambda^{*}$ and $\mu \leq \mu^{*},$ problem (1.5) has at last $\operatorname{cat}(\Omega)$ positive solutions for $4<p<6$ when $N=3$ or $2<p<2^{*}=\frac{2 N}{N-2}$ and $N \geq 4 .$ We recall that if $Y$ is a closed subset of a topological space $X$, the Lusternik-Schnirelman
category $\operatorname{cat}_X(Y)$ is the least number of closed and contractible sets in $X$ which cover $Y$. Hereafter, $\operatorname{cat}(X)$ denotes $\operatorname{cat}_X(X)$.
\par
Motivated by the above references, we intend to study the critical Schr\"odinger-Poisson system. Since in our case the problem is nonlocal, we need to be
careful in some estimates because some properties involving the energy functional that were used in the references mentioned above are not
clear that they work well in the present paper. We would like to point out that we focus on the case $N=3$ since such a problem stems from
physics, although the argument can be carried out for the space with different dimension.
\par
Now we state our main results as follows.
\begin{theorem}  Assume that $(V_1)-(V_3)$ are satisfied, then for $\lambda>0$ large enough, problem \eqref{maineq} has at least one positive ground state solution $u_{\lambda}$. Furthermore, for any sequence $\lambda_n\to +\infty,$
$\{u_{\lambda_n}\}$ has a subsequence converging to $u$ in $H^1(\mathbb{R}^3)$, where $u$ is a positive ground state solution of the limit problem
\begin{equation}\label{limiteq}
\begin{cases}
-\Delta u+\phi u=\mu |u|^{p-2}u+|u|^{4}u& \quad x\in \Omega,\\
-\Delta \phi=u^2&\quad x\in \Omega.\\
\end{cases}
\end{equation}
\end{theorem}

\begin{theorem}  Assume that $(V_1)-(V_3)$ are satisfied, then there exist $\lambda^*$, $\mu^*>0$ such that problem \eqref{maineq} has at least $\operatorname{cat}(\bar{\Omega})$ positive solutions for $\lambda\geq\lambda^*$ and $\mu\leq\mu^*$.
\end{theorem}

This paper is organized as follows. In the forthcoming section we collect some necessary preliminary Lemmas which will be used later. In section 3, we study the "limit" problem \eqref{limiteq}. In section 4, we are devoted to the existence of positive ground state solutions and their asymptotical behavior as $\lambda\to\infty$. In section 5, we discuss the relations among some mountain pass levels and complete the proof of multiple positive solutions.
\par
\vskip8pt
\textbf{Notation.}~In this paper we make use of the following notations.
\begin{itemize}
\item[$\bullet$] For any $R>0$ and for any $x\in\mathbb{R}^3$, $B_{R}(x)$ denotes the ball of radius $R$ centered at $x$.
%\item[$\bullet$]  $L^p(\mathbb{R}^3)$, $1\leq p<+\infty$ denotes the Lebesgue space with  the %norm $|u|_p=(\int_{\mathbb{R}^3}|u|^pdx)^{\frac{1}{p}}$.
\item[$\bullet$] The letters $C,C_i$ stand for positive constants (possibly different from line to line).
\item[$\bullet$]  "$\rightarrow$" for the strong convergence and "$\rightharpoonup$" for the weak convergence.
\item[$\bullet$]$|u|_q=(\int_{\mathbb{R}^3}|u|^qdx)^{\frac{1}{q}}$denotes the norm of $u$ in
$L^q(\mathbb{R}^3)$ for $2\leq q\leq 6$.
\end{itemize}

\section{Preliminaries}

In this section, we outline the variational framework for studying problem \eqref{maineq} and give the estimate of the energy level of the functional for which the $(PS)$ condition holds.  We first introduce some spaces. The standard norm of $H^1(\mathbb{R}^3)$ is given by
\begin{equation*}
\|u\|=\Big(\int_{\mathbb{R}^3}(|\nabla u|^2+u^2dx)\Big)^{1/2}.
\end{equation*}
Let
\begin{equation*}
E=\Big\{u\in H^1(\mathbb{R}^3)\ |\ \int_{\mathbb{R}^3}V(x)u^2dx<\infty\Big\},
\end{equation*}
be equipped with the inner product and norm
\begin{equation*}
(u,v)_\lambda=\int_{\mathbb{R}^3}(\nabla u\nabla v+\lambda V(x)uv)dx.
\end{equation*}
and
\begin{equation*}
\|u\|_\lambda=\Big(\int_{\mathbb{R}^3}(|\nabla u|^2+\lambda V(x)u^2)dx\Big)^{1/2}.
\end{equation*}
The conditions $(V_1)$ and $(V_2)$ yield that $E$ is a Hilbert space. Without loss of generality,
we assume $\lambda\geq 1$ throughout this paper. Then for any $u\in E$, there exists $\kappa>0$ such that
\begin{equation}\label{firstne}
\|u\|_\lambda\geq\kappa\|u\|.
\end{equation}
We also define Sobolev space $\mathcal{D}^{1,2}(\mathbb{R}^3)$ as the completion of $\mathcal{C}_0^\infty(\mathbb{R}^3)$ with respect to the norm
\begin{equation*}
\|u\|_{\mathcal{D}^{1,2}(\mathbb{R}^3)}=\Big(\int_{\mathbb{R}^3}(|\nabla u|^2dx\Big)^{\frac{1}{2}}.
\end{equation*}
It is clear that system \eqref{maineq} is the Euler-Lagrange equations of the functional $J:E\times \mathcal{D}^{1,2}(\mathbb{R}^3)\rightarrow\mathbb{R}$ defined by
\begin{equation*}
J(u,\phi)=\frac{1}{2}\|u\|^2_\lambda+\frac{1}{4}\int_{\mathbb{R}^3}\phi u^2dx-\frac{\mu}{p}\int_{\mathbb{R}^3}|u|^{p}dx-\frac{1}{6}\int_{\mathbb{R}^3}|u|^{6}dx.
\end{equation*}
It is easy to see that $J$ exhibits a strong indefiniteness, namely it is unbounded both from below and from above on infinitely dimensional subspaces. This indefiniteness can be removed using the reduction method described in \cite{Benci-Fortunato1998TMNA}. First of all, for a fixed $u\in H^1(\mathbb{R}^3)$ , there exists a unique $\phi_u\in \mathcal{D}^{1,2}(\mathbb{R}^3)$ which is the solution of
\[-\Delta\phi=u^2\ \ \text{in}\ \mathbb{R}^3.\]
We can write an integral expression for $\phi_u$ in the form
\begin{equation*}
\phi_u(x)=\frac{1}{4\pi}\int_{\mathbb{R}^3}\frac{ u^2(y)}{|x-y|}dy,\quad x\in \mathbb{R}^3,
\end{equation*}
which is called Riesz potential (see \cite{Landkof}).
Then the system \eqref{maineq} can be reduced to the first equation with $\phi$ represented by the solution of the  Poisson equation. This is the basic strategy of solving \eqref{maineq}. To be more precise about the solution $\phi$ of the  Poisson equation, we collect some useful lemmas.

\begin{lemma}\cite{zhaozhao2009}\label{zhao1}
For any $u\in H^1(\mathbb{R}^3)$, we have
\begin{itemize}
\item[$(i)$] $\phi_u\geq0$;
\item[$(ii)$]$\phi_u:H^1(\mathbb{R}^3)\rightarrow \mathcal{D}^{1,2}(\mathbb{R}^3)$ is continuous and maps bounded sets into bounded sets;
\item[$(iii)$]$\|\phi_u\|^2_{D^{1,2}(\mathbb{R}^3)}=\int_{\mathbb{R}^3}\phi_u u^2dx\leq C|u|_{\frac{12}{5}}^4$;
\item[$(iv)$]If $u_n\rightharpoonup u$ in $H^1(\mathbb{R}^3)$, then $\phi_{u_n} \rightharpoonup \phi_u$ in $\mathcal{D}^{1,2}(\mathbb{R}^3)$;
\item[$(v)$]If $u_n\rightarrow u$ in $H^1(\mathbb{R}^3)$, then $\phi_{u_n} \rightarrow \phi_u$ in $\mathcal{D}^{1,2}(\mathbb{R}^3)$ and $\int_{\mathbb{R}^3}\phi_{u_n}u_n^2dx\rightarrow \int_{\mathbb{R}^3}\phi_u u^2dx$.
\end{itemize}
\end{lemma}

Define $F:H^1(\mathbb{R}^3)\rightarrow\mathbb{R}$ by
\[F(u)=\int_{\mathbb{R}^3}\phi_u u^2dx.\]
It is apparent that $F(u(\cdot+y))=F(u)$ for any $y\in\mathbb{R}^3$, $u\in H^1(\mathbb{R}^3)$ and $F$ is weakly lower semi-continuous in $H^1(\mathbb{R}^3)$. Moreover, similarly to the well-know Brezis-Lieb Lemma (\cite{Brezis-Lieb1983PAMS}), we have the next Lemma.

\begin{lemma}\cite{zhaozhao2009}\label{zhao2} Let $u_n\rightharpoonup u$ in $H^1(\mathbb{R}^3)$ and $u_n\rightarrow u$ a.e.in $\mathbb{R}^3$. Then
\begin{itemize}
\item[$(i)$]$F(u_n-u)=F(u_n)-F(u)+o(1)$;
\item[$(ii)$]$F^\prime(u_n-u)=F^\prime(u_n)-F^\prime(u)+o(1)$, in $(H^1(\mathbb{R}^3))^{-1}$.
\end{itemize}
\end{lemma}

Putting $\phi=\phi_u$ into the first equation of \eqref{maineq}, we obtain a nonlocal semilinear elliptic equation
\[-\Delta u+\lambda V(x)u+\phi_u u=\mu |u|^{p-2}u+|u|^{4}u~~\text{in}\ \mathbb{R}^3.\]
The corresponding functional is
\begin{equation*}
I_{\lambda,\mu}(u)=\frac{1}{2}\|u\|^2_\lambda+\frac{1}{4}F(u)-\frac{\mu}{p}\int_{\mathbb{R}^3}|u|^{p}dx-\frac{1}{6}\int_{\mathbb{R}^3}|u|^{6}dx.
\end{equation*}
It is easy to check that $I_{\lambda,\mu}$ is well defined on $E$ and $I_{\lambda,\mu} \in C^1(E,\mathbb{R}).$ Then we can define
\begin{equation*}
\mathcal{N}_{\lambda,\mu}=\{u\in E\setminus \{0\}\ |\  \langle I'_{\lambda,\mu}(u), u\rangle =0\}.
\end{equation*}

\begin{lemma}\label{le:lowboundedinN}There exists $\sigma>0$ which is independent of
 $\lambda$ such that
\begin{equation*}
\|u\|_\lambda>\sigma \quad\text{and}\quad I_{\lambda,\mu}(u)\geq\frac{p-2}{2p}\sigma^2, \quad\text{for all}\ u\in \mathcal{N}_{\lambda,\mu}.
\end{equation*}
\end{lemma}
\noindent{\bf Proof:} For any $u\in N_{\lambda,\mu}$, without loss of generality, we assume that $\|u\|_\lambda\leq 1$.  From \eqref{firstne},  we have
\begin{equation*}
\begin{aligned}
\displaystyle\|u\|^2_\lambda+F(u)&=\displaystyle\mu\int_{\mathbb{R}^3}
|u|^pdx+\int_{\mathbb{R}^3}|u|^6dx.\\
&\leq C(\|u\|^p_\lambda+\|u\|^6_\lambda).
\end{aligned}
\end{equation*}
Then the first desired result follows from $\|u\|_\lambda\leq 1$.
\par
On the other hand,  we have
\begin{equation*}
\begin{aligned}
I_{\lambda,\mu}(u)&=\displaystyle\frac{1}{2}\|u\|^2_\lambda+\frac{1}{4}F(u)-\frac{\mu}{p}\int_{\mathbb{R}^3}|u|^{p}dx-\frac{1}{6}\int_{\mathbb{R}^3}|u|^{6}dx\\
&\geq\displaystyle\frac{1}{2}\|u\|^2_\lambda+\frac{1}{4}F(u)-\frac{\mu}{p}\int_{\mathbb{R}^3}|u|^{p}dx-\frac{1}{p}\int_{\mathbb{R}^3}|u|^{6}dx\\
&=\displaystyle\frac{1}{2}\|u\|^2_\lambda+\frac{1}{4}F(u)-\frac{1}{p}(\|u\|^2_\lambda+F(u))\\
&\geq\displaystyle(\frac{1}{2}-\frac{1}{p})\|u\|^2_\lambda\\
&\geq\displaystyle\frac{p-2}{2p}\sigma^2.
\end{aligned}
\end{equation*}
\hfill{$\Box$}

\begin{lemma}\label{le:nehari}
For any $u\in E\setminus\{0\}$, there exists a unique $t(u)>0$ such that $t(u)u\in\mathcal{N}_{\lambda,\mu}$ and \[I_{\lambda,\mu}(t(u)u)=\displaystyle\max_{t\geq0}I_{\lambda,\mu}(tu).\]
\end{lemma}
{\bf Proof:} For any $u\in E\setminus\left\{0\right \}$, define $g(t)=I_{\lambda,\mu}(tu),\ t\in[0,+\infty).$ Then
\[g(t)=\frac{t^2}{2}\|u\|^2_\lambda+\frac{t^4}{4}F(u)-\frac{t^p\mu}{p}\int_{\mathbb{R}^3}|u|^{p}dx-\frac{t^6}{6}\int_{\mathbb{R}^3}|u|^{6}dx.\]
It is easy to see that $g(t)>0$ for $t>0$ small and $g(t)<0$ for $t>0$ large enough, so there exists $t_0>0$ such that
\[g'(t_0)=0\quad \hbox{and}\quad g(t_0)=\max_{t\geq 0}g(t)=\max_{t\geq 0}I_{\lambda,\mu}(tu).\]
It follows from $g'(t_0)=0$ that $t_0u\in \mathcal{N}_{\lambda,\mu}$.

If there exist $0<t_1<t_2$ such that $t_1u\in \mathcal{N}_{\lambda,\mu}$ and $t_2u\in \mathcal{N}_{\lambda,\mu}$. Then
\[\frac{1}{t_1^2}\|u\|^2_{\lambda}+F(u)=t_1^{p-4}\mu\int_{\mathbb{R}^3}|u|^{p}dx+t_1^2\int_{\mathbb{R}^3}|u|^{6}dx\]
and
\[\frac{1}{t_2^2}\|u\|^2_{\lambda}+F(u)=t_2^{p-4}\mu\int_{\mathbb{R}^3}|u|^{p}dx+t_2^2\int_{\mathbb{R}^3}|u|^{6}dx.\]
It follows that
$$(\frac{1}{t_1^2}-\frac{1}{t_2^2})\|u\|^2_{\lambda}=(t_1^{p-4}-t_2^{p-4})\mu\int_{\mathbb{R}^3}|u|^{p}dx+(t_1^2-t_2^2)\int_{\mathbb{R}^3}|u|^{6}dx,$$
which is a contradiction.
\hfill{$\Box$}

%引理5

\begin{lemma}\label{le:c-equal}
For any $\lambda\geq 1$ and $\mu>0$, let \[c_{\lambda,\mu}=\inf_{u\in \mathcal{N}_{\lambda,\mu}}I_{\lambda,\mu}(u),\quad c_{\lambda,\mu}^*=\inf_{u\in E\setminus\{0\}}\max_{t\geq0}I_{\lambda,\mu}(tv),
\quad c_{\lambda,\mu}^{**}=\inf_{\gamma\in \Gamma}\sup_{t\in [0,1]}I_{\lambda,\mu}(\gamma(t)),\]
where
\[\Gamma=\{\gamma(t)\in C([0,1],E)\ |\ \gamma(0)=0,\  I_{\lambda,\mu}(\gamma(1))<0\}.\] Then, $c_{\lambda,\mu}=c_{\lambda,\mu}^*=c_{\lambda,\mu}^{**}.$
\end{lemma}
{\bf Proof:} We  divide the proof into three steps.

{\bf Step1}. $c_{\lambda,\mu}^*=c_{\lambda,\mu}.$ By Lemma  \ref{le:nehari}, we have \[c_{\lambda,\mu}^*=\displaystyle \inf_{u\in E\setminus\{0\}}\max_{t\geq0}I_{\lambda,\mu}(tu)
=\inf_{u\in E\setminus\{0\}}I_{\lambda,\mu}(t(u)u)=\inf_{u\in \mathcal{N}_{\lambda,\mu}}I_{\lambda,\mu}(u)=c_{\lambda,\mu}.\]

{\bf Step2}. $c_{\lambda,\mu}^*\geq c_{\lambda,\mu}^{**}.$
From Lemma \ref{le:nehari}, for any $u\in E\setminus\{0\}$, there exists $T$ large enough, such that $I_{\lambda,\mu}(Tu)<0.$ Define $\gamma(t)=tTu$, $t\in[0,1]$. Then we have
$\gamma(t)\in \Gamma$ and,  therefore,
\[c_{\lambda,\mu}^{**}=\inf_{\gamma\in \Gamma}\sup_{t\in [0,1]}I_{\lambda,\mu}(\gamma(t))\leq \sup_{t\in [0,1]}I_{\lambda,\mu}(\gamma(t))\leq \max_{t\geq0}I_{\lambda,\mu}(tu).\]
It follows that $c_{\lambda,\mu}^*\geq c_{\lambda,\mu}^{**}$.

{\bf Step3}.  $ c_{\lambda,\mu}^{**}\geq c_{\lambda,\mu}.$
For any $u\in E\setminus \{0\}$ with
$\|u\|_\lambda$ small, we know
\begin{equation}\label{eq1}
\|u\|^2_\lambda+F(u)>\mu\int_{\mathbb{R}^3}|u|^pdx+\int_{\mathbb{R}^3}|u|^6dx.
\end{equation}
We claim that every $\gamma(t)\in \Gamma$ has to cross $N_{\lambda,\mu}$. Otherwise, by the continuity of $\gamma(t)$, \eqref{eq1} still holds, when $u$ is replaced by $\gamma(1)$.  Then we can obtain
\begin{equation*}
\begin{aligned}
I_{\lambda,\mu}(\gamma(1))&=\displaystyle\frac{1}{2}\|\gamma(1)\|^2_\lambda+\frac{1}{4}F(\gamma(1))-\frac{\mu}{p}\int_{\mathbb{R}^3}|\gamma(1)|^{p}dx-\frac{1}{6}\int_{\mathbb{R}^3}|\gamma(1)|^{6}dx\\
&\geq\displaystyle\frac{1}{2}\|\gamma(1)\|^2_\lambda+\frac{1}{4}F(\gamma(1))-\frac{\mu}{p}\int_{\mathbb{R}^3}|\gamma(1)|^{p}dx-\frac{1}{p}\int_{\mathbb{R}^3}|\gamma(1)|^{6}dx\\
&\geq\displaystyle\frac{1}{2}\|\gamma(1)\|^2_\lambda+\frac{1}{4}F(\gamma(1))-\frac{1}{p}(\|\gamma(1)\|^2_\lambda+F(\gamma(1)))\\
&\geq\displaystyle(\frac{1}{2}-\frac{1}{p})\|\gamma(1)\|^2_\lambda\\
&> 0,
\end{aligned}
\end{equation*}
which contradicts the definition of $\gamma(1)$. It follows from the claim that  $ c_{\lambda,\mu}^{**}\geq c_{\lambda,\mu}$.
\hfill{$\Box$}
\par
One can easily check that the functional $I_{\lambda,\mu}$ satisfies the mountain-pass geometry, that is the following lemma holds.
\begin{lemma}\label{le:moutain-pass geometry}
$I_{\lambda,\mu}$ has the mountain geometry structure.
\begin{itemize}
\item[(1)] There exist $a_0,r_0>0$ independent of $\lambda$, such that $I_{\lambda,\mu}(u)\geq a_0$, for all $u\in E$ with $\|u\|_\lambda=r_0.$
\item[(2)] For any $u\in E\setminus\{0\}$ , $\lim_{t\to \infty}I_{\lambda,\mu}(tu)=-\infty.$
\end{itemize}
\end{lemma}

\begin{lemma}\label{le:c-level}
For any $\lambda\geq 1$, $\mu>0$, we have $c_{\lambda,\mu}<\displaystyle\frac{1}{3} S^{3/2}$.
\end{lemma}
{\bf Proof:}
Given $\epsilon>0,$ we consider the function
\[U_\epsilon(x)=\displaystyle\frac{(3\epsilon)^{\frac{1}{4}}}{(\epsilon+|x|^2)^{\frac{1}{2}}}\]
which satisfies
\[\int_{\mathbb{R}^3}|\nabla U_\epsilon|^2dx=\int_{\mathbb{R}^3}| U_\epsilon|^6 dx=S^{\frac{3}{2}}.\]
Let $\varphi\in C_0^\infty(\mathbb{R}^3,[0,1])$ be such that $\varphi=1$ if $|x|<1$ and  $\varphi=0$ if $|x|\geq2$.
Consider the following function
$v_\epsilon(x):=\varphi(x) U_\epsilon(x).$
By direct computation, we get
\begin{equation}\label{est1}
|\nabla v_\epsilon|^{2}_2=S^{\frac{3}{2}}+O(\epsilon^{\frac{1}{2}}),
\end{equation}
\begin{equation}\label{est2}
|v_\epsilon|^{2}_6=S^{\frac{1}{2}}+O(\epsilon).
\end{equation}
\begin{equation}\label{est3}
| v_\epsilon|^{r}_r\sim\left\{
\begin{array}{ll} O(\epsilon^{\frac{r}{4}}), &\quad r\in [2,3),\\
O(\epsilon^{\frac{r}{4}}|ln\epsilon|), &\quad r=3,\\
O(\epsilon^{\frac{6-r}{4}}), &\quad r\in (3,6).
\end{array}\right.
\end{equation}
By Lemma \ref{le:nehari}, there exists $t_\epsilon>0$ such that
$t_\epsilon v_\epsilon\in \mathcal{N}_{\lambda,\mu}$ and  \[I_{\lambda,\mu}(t_\epsilon v_\epsilon)=\displaystyle\max_{t\geq0}I_{\lambda,\mu}(t v_\epsilon).\]
From Lemma \ref{le:lowboundedinN}, we have $\|t_\epsilon v_\epsilon\|_\lambda\geq \sigma$. Then it follows from \eqref{est1} and \eqref{est3} that $t_\epsilon$ has a positive lower bound. On the other hand, in view of $t_\epsilon v_\epsilon\in \mathcal{N}_{\lambda,\mu}$, we can obtain
\begin{equation*}
\begin{aligned}
\displaystyle t_\epsilon^6\int_{\mathbb{R}^3}|v_\epsilon|^6dx&\leq\displaystyle t_\epsilon^2 \|v_\epsilon\|^2_\lambda+t_\epsilon^4 F(v_\epsilon)\\
&\leq  t_\epsilon^2 \|v_\epsilon\|^2_\lambda+t_\epsilon^4 |v_\epsilon|^4_{\frac{12}{5}}.
\end{aligned}
\end{equation*}
Noting \eqref{est1}, \eqref{est2} and \eqref{est3}, we can know that $t_\epsilon$ is bounded from above.  Therefore,
\begin{equation*}
\begin{aligned}
c_{\lambda,\mu}\leq I_{\lambda,\mu}(t_\epsilon v_\epsilon)&=\displaystyle\frac{t_\epsilon^2}{2}\|v_\epsilon\|^2_\lambda+\frac{t_\epsilon^4}{4}F(v_\epsilon)-\frac{t_\epsilon^p}{p}\int_{\mathbb{R}^3}|v_\epsilon|^{p}dx-\frac{t_\epsilon^6}{6}\int_{\mathbb{R}^3}|v_\epsilon|^{6}dx\\
&\leq \displaystyle\frac{t_\epsilon^2}{2}(S^{\frac{3}{2}}+C_1\epsilon^{\frac{1}{2}}+C_2\epsilon^{\frac{1}{2}})+C_3\epsilon-C_4\epsilon^{\frac{6-p}{4}}-\frac{t_\epsilon^6}{6}(S^{\frac{1}{2}}-C_5\epsilon)^3\\
&\leq \displaystyle (\frac{t_\epsilon^2}{2}-\frac{t_\epsilon^6}{6})S^{\frac{3}{2}}+C_6\epsilon^{\frac{1}{2}}-C_4\epsilon^{\frac{6-p}{4}}\\
&\leq \displaystyle \frac{1}{3}S^{\frac{3}{2}}+C_6\epsilon^{\frac{1}{2}}-C_4\epsilon^{\frac{6-p}{4}}\\
&< \displaystyle \frac{1}{3}S^{\frac{3}{2}}.
\end{aligned}
\end{equation*}
\hfill{$\Box$}

%引理8

\begin{lemma}\label{le:psbounded} Any of the  $(PS)_{c_{\lambda,\mu}}$ sequence $\left\{u_n\right\}$
for $I_{\lambda,\mu}$ is bounded, and
\[\limsup_{n\to \infty} \|u_n\|_{\lambda}\leq
\sqrt{\frac{2p}{p-2}}c_{\lambda,\mu}.\]
\end{lemma}
{\bf Proof:} Suppose that $\left\{u_n\right\}$ is a $(PS)_{c_{\lambda,\mu}}$ sequence of
$I_{\lambda,\mu}$,  we have \[I_{\lambda,\mu}(u_n)\to _{c_{\lambda,\mu}},\quad I'_{\lambda,\mu}(u_n)\to 0 .\]
Thus
\begin{equation*}
\begin{aligned}
c_{\lambda,\mu}+o(1)+o(1)\|u_n\|_\lambda
&=I_{\lambda,\mu}(u_n)-\displaystyle \frac{1}{p}\langle I'_{\lambda,\mu}(u_n),
v_n\rangle \\
&=\displaystyle(\frac{1}{2}-\frac{1}{p})\|u_n\|^2_\lambda+(\frac{1}{4}-\frac{1}{p})F(u_n)\\
&\quad\displaystyle+(\frac{1}{p}-\frac{1}{6})\int_{\mathbb{R}^3}
u_n^6dx.
\end{aligned}
\end{equation*}
It follows that
\[(\frac{1}{2}-\frac{1}{p})\|u_n\|^2_\lambda\leq c_{\lambda,\mu}+o(1)+o(1)\|u_n\|_\lambda.\]
Then $\left\{u_n\right\}$ is bounded in $E$, and
\[\limsup_{n\to \infty} ||u_n||_\lambda\leq
\sqrt{\frac{2p}{p-2}}c_{\lambda,\mu}.\]
\hfill{$\Box$}

%引理9

\begin{lemma}\label{prop:strong} Let $M>0$ be a constant, and as $n\to\infty$, $\lambda_n\rightarrow\infty$. If $\|u_n\|_{\lambda_n}\leq M$ and $u_n\rightharpoonup 0$ in $H^1(\mathbb{R}^3)$, then $\displaystyle\lim_{n\to \infty}\int_{\mathbb{R}^3}|u_n|^{q}dx=0$, for $2<q<6.$
\end{lemma}
{\bf Proof:} For any $R>0$,  let
\[
A(R):=\left\{x\in\mathbb{R}^3| |x|\geq R, \ V(x)\geq M_0\right\},
\]
\[
B(R):=\left\{x\in\mathbb{R}^3| |x|\geq R, \ V(x)\leq M_0\right\}.
\]
It is easy to see that
\begin{equation}\label{foutsmall}
\begin{array}{ll}
\displaystyle\int_{A(R)}u^2_ndx &\leq \displaystyle\frac{1}{\lambda_nM_0}
\int_{A(R)}\lambda_n V(x)u^2_ndx\\
&\leq \displaystyle\frac{1}{\lambda_nM_0}
\int_{A(R)}(|\nabla u_n|^2+\lambda_n V(x)u^2_n)dx\\
&=\displaystyle\frac{1}{\lambda_nM_0}\|u_n\|^2_{\lambda_n} \\
&\leq\displaystyle\frac{1}{\lambda_nM_0}M^2
\end{array}
\end{equation}
and
\begin{equation}\label{foutsmall2}
\begin{array}{ll}
\displaystyle\int_{B(R)}u^2_ndx&\leq
\displaystyle [m(B(R))]^{\frac{2}{3}}\Big(\int_{\mathbb{R}^3}u^6_n dx\Big)^{\frac{1}{3}}\\
&\leq\displaystyle [m(B(R))]^{\frac{2}{3}}S^{-1}\int_{\mathbb{R}^3}|\nabla u_n|^2dx\\
&\leq\displaystyle M^2S^{-1}[m(B(R))]^{\frac{2}{3}}.
\end{array}
\end{equation}
By using interpolation inequality, there exists $\theta\in(0,1)$ such that
\begin{equation*}
\begin{array}{ll}
\displaystyle\int_{B_R^c}|u_n|^{q}dx&\leq
\displaystyle(\int_{B_R^c}|u_n|^2dx)^{\frac{q\theta}{2}}(\int_{B_R^c}|
u_n|^{6} dx)^\frac{q(1-\theta)}{6}\\
&\leq S^{-\frac{q(1-\theta)}{2}}\displaystyle(\int_{B_R^c}|u_n|^2dx)^{\frac{q\theta}{2}}(\int_{\mathbb{R}^N}|\nabla u_n|^2 dx)^\frac{q(1-\theta)}{2}\\
&\leq \displaystyle S^{-\frac{q(1-\theta)}{2}}M^{q(1-\theta)}(\int_{A(R)}|u_n|^2dx+\int_{B(R)}|
u_n|^2dx)^{\frac{q\theta}{2}}.
\end{array}
\end{equation*}
Then from \eqref{foutsmall} and \eqref{foutsmall2}, there exists $C>0$ such that
\begin{equation}\label{foutsmall3}
\int_{B_R^c}|u_n|^{q}dx\leq C\Big(\frac{1}{\lambda_n}+[m(B(R))]^{\frac{2}{3}}\Big).
\end{equation}
In view of $(V_2)$ and \eqref{foutsmall3}, for any
$\epsilon>0,$ there exist $N_\epsilon, R_\epsilon>0$ such that
\begin{equation}\label{foutsmall4}
\int_{B_R^c}|u_n|^{q}dx\leq \epsilon,
\end{equation}
when $n\geq N_\epsilon $ and $R\geq R_\epsilon$. Fix $R_1>R_\epsilon$. Then, for $n\geq N_\epsilon $, we have
\begin{equation*}
\begin{array}{ll}
\displaystyle\int_{\mathbb{R}^3}|u_n|^{q}dx&=\displaystyle\int_{B_{R_1}}|u_n|^{q}dx+\int_{B^c_{R_1}}|u_n|^{q}dx\\
&\leq \displaystyle\int_{B_{R_1}}|u_n|^{q}dx+\epsilon.
\end{array}
\end{equation*}
Noting  $u_n\rightharpoonup 0$ in $H^1(\mathbb{R}^3)$, it follows from the above inequality that
\[\displaystyle\lim_{n\to \infty}\int_{\mathbb{R}^3}|u_n|^{q}dx=0.\]
\hfill{$\Box$}

To find a new upper bound for $c_{\lambda,\mu}$, we need to investigate the "limit" problem \eqref{limiteq}. Firstly,  denote the standard norm of $H_0^1(\Omega)$ by
\begin{equation*}
\|u\|_0:=\Big(\int_{\Omega}|\nabla u|^2dx)\Big)^{1/2}.
\end{equation*}
The associated energy functional is
\begin{equation*}
I_{\mu}(u)=\frac{1}{2}\|u\|^2_0+\frac{1}{4}F_0(u)-\frac{\mu}{p}\int_{\Omega}|u|^{p}dx-\frac{1}{6}\int_{\Omega}|u|^{6}dx,
\end{equation*}
where $F_0(u)=\frac{1}{4\pi}\iint_{\Omega\times\Omega}\frac{u^2(x)u^2(y)}{|x-y|}dxdy$.
It is easy to check that $I_{\mu}$ is well defined on $H_0^1(\Omega)$ and $I_{\mu} \in C^1(H_0^1(\Omega),\mathbb{R})$. Then we can define
\begin{equation*}
\mathcal{N}_{\mu}=\{u\in H_0^1(\Omega)\setminus \{0\}\ |\  \langle I'_{\mu}(u), u\rangle =0\}
\end{equation*}
and
\begin{equation*}
c_{\mu}=\inf_{u\in \mathcal{N}_{\mu}}I_{\mu}(u)
\end{equation*}

\begin{remark}\label{ro}
For $c_{\mu}$, $I_{\mu}$ and $\mathcal{N}_{\mu}$, there are silimar results obtained from
lemma \ref{le:lowboundedinN} to Lemma \ref{le:c-level}. By Mountain pass Theorem, we can see that there exists $u\in H_0^1(\Omega)$ such that
$I_{\mu}(u)=c_{\mu}$ and $I'_{\mu}(u)=0$.
\end{remark}

\begin{lemma}\label{cbounded} For $\lambda\geq 1$, $\mu>0$, we have $c_{\lambda,\mu}\leq c_{\mu}$.
\end{lemma}
{\bf Proof:} For $u\in \mathcal{N}_\mu$, we have
\begin{equation*}
\int_{\Omega}|\nabla u|^2dx+\int_{\Omega}\phi_uu^2dx
=\mu\int_{\Omega}|u|^{p-2}udx+\int_{\Omega}|u|^{4}udx.
\end{equation*}
By $V(x)=0$ in $\Omega$ and $u=0$ in $\mathbb{R}^N\setminus\Omega,$ the above equality can be written as
\begin{equation*}
\int_{\mathbb{R}^3}(|\nabla u|^2+\lambda V(x)u^2)dx+\int_{\mathbb{R}^3}\phi_uu^2dx
=\mu\int_{\mathbb{R}^3}|u|^{p-2}udx+\int_{\mathbb{R}^3}|u|^{4}udx.
\end{equation*}
Thus $ u\in \mathcal{N}_{\lambda,\mu}$.

On the other hand, we have
\begin{equation*}
\begin{aligned}
I_\mu(u)&=\displaystyle\frac{1}{2}\int_{\Omega}|\nabla
u|^2dx+\frac{1}{4}\int_{\Omega}\phi_uu^2dx-\frac{\mu}{p}\int_{\Omega}|u|^pdx-\frac{1}{6}\int_{\Omega}|u|^6dx\\
&=\displaystyle\frac{1}{2}\int_{\mathbb{R}^3}(|\nabla u|^2+\lambda V(x)u^2)dx+\frac{1}{4}\int_{\mathbb{R}^3}\phi_uu^2dx-\frac{\mu}{p}\int_{\mathbb{R}^3}|u|^pdx-\frac{1}{6}\int_{\Omega}|u|^6dx\\
&=I_{\lambda,\mu}(u).
\end{aligned}
\end{equation*}
Therefore, $c_{\lambda,\mu}\leq c_\mu$.
\hfill{$\Box$}

\begin{corollary}\label{r1}
It follows from Remark \ref{ro} that there exists $\tau>0,$ such that $c_{\lambda,\mu}<\displaystyle\frac{1}{3}S^{\frac{3}{2}}-\tau,$ for $\lambda\geq 1$, $\mu>0$.
\end{corollary}

\begin{lemma}\label{le:pscondition}
There is $\lambda^*>0$ such that $I_{\lambda,\mu}$ verifies $(PS)_{d_\lambda}$ condition
for any $d_\lambda\in(0, \frac{1}{3}S^{\frac{3}{2}}-\tau)$ for all $\lambda\geq\lambda^*$.
\end{lemma}
{\bf Proof:} Let $\left\{u_n\right\}$ be a $(PS)_{d_\lambda}$ sequence of $I_{\lambda,\mu}$, that is
\[I_{\lambda,\mu}(u_n)\to d_\lambda,\quad I'_{\lambda,\mu}(u_n)\to 0.\]
By Lemma \ref{le:psbounded},  $\left\{u_n\right\}$ is bounded in $E.$ Thus, there is $u\in E$ such taht
\begin{equation*}
\begin{aligned}
&u_n\rightharpoonup u\quad\hbox{in}\; E,\\
&u_n\to u\quad\hbox{a.e. in}\  \mathbb{R}^3,\\
&u_n\rightharpoonup u
\quad\hbox{in}\  L^{q}(\mathbb{R}^3), \ \hbox{for }\ 2\leq q \leq 6.
\end{aligned}
\end{equation*}
Define $v_n=u_n-u$. From Lemma \ref{zhao2}, we know $\{v_n\}$ is also a $PS$ sequence for $I_{\lambda,\mu}$ and
\[I_{\lambda,\mu}(v_n)\to d_\lambda-I_{\lambda,\mu}(u),\quad I'_{\lambda,\mu}(v_n)\to 0.\]
By using $I'_{\lambda,\mu}(v_n)\to 0$ as $n\rightarrow\infty$, we have
\begin{equation}\label{lim20}
o_n(1)=\langle I'_{\lambda,\mu}(v_n),v_n\rangle
=\|v_n\|^2_{\lambda}+F(v_n)-\mu\int_{\mathbb{R}^3}|v_n|^{p}dx-\int_{\mathbb{R}^3}|v_n|^{6}dx.
\end{equation}
From the boundness of $\left\{v_n\right\}$,  we can assume
\begin{equation}\label{lim14}
\begin{array}{ll}
\|v_n\|^2_{\lambda}\to L^{(1)}_\lambda,
\end{array}
\end{equation}
\begin{equation}\label{lim24}
\begin{array}{ll}
F(v_n)\to L^{(2)}_\lambda,
\end{array}
\end{equation}
and
\begin{equation}\label{lim34}
\mu\int_{\mathbb{R}^3}|v_n|^{p}+\int_{\mathbb{R}^3}|v_n|^6\to L^{(3)}_\lambda.
\end{equation}
Then, we have
\begin{equation}\label{lim54}
L^{(1)}_\lambda+L^{(2)}_\lambda= L^{(3)}_\lambda.
\end{equation}
If $L^{(1)}_\lambda=0$,  we have $v_n\rightarrow 0$ in $E$ which means
$u_n\rightarrow u$ in $E.$ Next, we prove that when $\lambda$ is large enough, $L^{(1)}_\lambda=0.$ If not,  we assume
\[\mu\int_{\mathbb{R}^3}|v_n|^{p}dx\to A_\lambda \quad\hbox{and}\quad \int_{\mathbb{R}^3}|v_n|^6\to B_\lambda.\]
Then, we can get
\begin{equation}\label{lim64}
L^{(1)}_\lambda+L^{(2)}_\lambda= L^{(3)}_\lambda= A_\lambda+B_\lambda.
\end{equation}
Arguing as in the proof of Lemma \ref{prop:strong}, we see that $A_\lambda=o_\lambda(1)$ and $L^{(2)}_\lambda=o_\lambda(1)$ as $\lambda\rightarrow\infty$. Therefore,
\begin{equation}\label{limit84}
L^{(1)}_\lambda=B_\lambda+o_\lambda(1).
\end{equation}
From \eqref{lim20}, there exists $C>0$ such that
\[
\|v_n\|^2_{\lambda}\leq C\Big(\|v_n\|_{\lambda}^{p}+\|v_n\|_{\lambda}^{6}+o_n(1)\Big).
\]
If $L^{(1)}_\lambda\leq1,$ letting $n\to \infty$ in the above inequality,
\begin{equation*}
L^{(1)}_\lambda\leq C\Big((L^{(1)}_\lambda)^{\frac{p}{2}}+(L^{(1)}_\lambda)^{3}\Big).
\end{equation*}
Thus, there exists $C>0$ independent of $\lambda$ such that
\begin{equation}\label{lim74}
L^{(1)}_\lambda\geq C.
\end{equation}
By the definition of $S$, we can get
\[S(\int_{\mathbb{R}^3}|v_n|^{6}dx)^{\frac{1}{3}}\leq\|v_n\|^2_\lambda.\]
Letting $n\to \infty$, we have
\begin{equation}\label{4equality2}
SB_\lambda^{\frac{1}{3}}\leq L^{(1)}_\lambda.
\end{equation}
Using \eqref{limit84}, \eqref{lim74} and \eqref{4equality2}, we obtain  \[\liminf_{\lambda\rightarrow\infty}L^{(1)}_\lambda\geq S^{\frac{3}{2}}.\]
From $I_{\lambda,\mu}(v_n)\to d_\lambda-I_{\lambda,\mu}(u)$, we have
\[
d_\lambda-I_{\lambda,\mu}(u)=\frac{1}{2}\|v_n\|^2_{\lambda}+\frac{1}{4}F(v_n)-\frac{\mu}{p}\int_{\mathbb{R}^3}|v_n|^{p}dx-\frac{1}{6}\int_{\mathbb{R}^3}|v_n|^{6}dx+o_n(1).
\]
Letting $n\rightarrow\infty$, we get
\[d_\lambda-I_{\lambda,\mu}(u)=\frac{1}{2}L^{(1)}_\lambda+\frac{1}{4}L^{(2)}_\lambda-\frac{1}{p}A_\lambda-\frac{1}{6}B_\lambda.\]
Noting $I'_{\lambda,\mu}(u)=0$, by Lemma \ref{le:lowboundedinN}, we have $I_{\lambda,\mu}(u)\geq 0$. Therefore,
\begin{equation*}
\begin{aligned}
\displaystyle\frac{1}{3}S^{\frac{3}{2}}-\tau&\geq\displaystyle\liminf_{\lambda\rightarrow\infty}\Big(\frac{1}{2}L^{(1)}_\lambda+\frac{1}{4}L^{(2)}_\lambda-\frac{1}{p}A_\lambda-\frac{1}{6}B_\lambda\Big)\\
&\geq\displaystyle(\frac{1}{2}-\frac{1}{6})S^{\frac{3}{2}}\\
&=\displaystyle\frac{1}{3}S^{\frac{3}{2}}
\end{aligned}
\end{equation*}
which is a contradiction.
\hfill{$\Box$}

\begin{corollary}\label{le:nepscondition}
There is $\lambda^*>0$ such that $I_{\lambda,\mu}$ verifies $(PS)_{d_\lambda}$ condition
on $\mathcal{N}_{\lambda,\mu}$ for any $d_\lambda\in(0, \frac{1}{3}S^{\frac{3}{2}}-\tau)$ for all $\lambda\geq\lambda^*$.
\end{corollary}

\begin{lemma}\label{le:same} If $u$ is a critical point of $I_{\lambda,\mu}$ on $\mathcal{N}_{\lambda,\mu}$, then $u$ is a critical point of $I_{\lambda,\mu}$ in $E$.
\end{lemma}
{\bf Proof:} Since $u$ is a critical point of $I_{\lambda,\mu}$ on $\mathcal{N}_{\lambda,\mu}$,  there exists $\theta\in \mathbb{R}$ such that
\[I'_{\lambda,\mu}(u)=\theta J'_{\lambda,\mu}(u)\]
where $J_{\lambda,\mu}(u)=\langle I'_{\lambda,\mu}(u), u\rangle$.

It follows from  $u\in N_{\lambda,\mu}$ that
\begin{equation*}
\|u\|^2_\lambda+F(u)=\mu\int_{\mathbb{R}^3}|u|^{p}dx+\int_{\mathbb{R}^3}|u|^{6}dx.
\end{equation*}
Then
\begin{equation*}
\begin{aligned}
\langle J'_{\lambda,\mu}(u), u\rangle&=\displaystyle 2\|u\|^2_{\lambda}+4F(u)- p\mu\int_{\mathbb{R}^3}|u|^{p}dx-6\int_{\mathbb{R}^3}|u|^{6}dx\\
&=\displaystyle -2\|u\|^2_{\lambda}+(4-p)\mu\int_{\mathbb{R}^3}|u|^{p}dx-2\int_{\mathbb{R}^3}|u|^{6}dx<0.
\end{aligned}
\end{equation*}
From $0=\langle I'_{\lambda,\mu}(u), u\rangle=\theta\langle J'_{\lambda,\mu}(u), u\rangle$ and the above inequality, we have $I'_{\lambda,\mu}(u)=0$.
\hfill{$\Box$}

\section{The asymptotical behavior of solutions}

%引理3.1

\begin{lemma}\label{le:clambdaachieved}
When $\lambda$ is large enough, problem \eqref{maineq} has at least one positive ground state solution.
\end{lemma}
{\bf Proof:} By Lemma \ref{le:c-equal} and Lemma \ref{le:moutain-pass geometry}, there exits a sequence $\left\{u_n\right\}$ which is a $(PS)_{c_{\lambda,\mu}}$ sequence of $I_{\lambda,\mu}$. Hence, by Corollary \ref{r1} and Lemma \ref{le:pscondition}, when $\lambda$ is large enough, (up to a subsequence) we have $u_n\rightarrow u$ in $E$. Therefore, we have $I_{\lambda,\mu}(u)=c_{\lambda,\mu}$ and $I'_{\lambda,\mu}(u)=0.$
Moreover, $|u|\in N_{\lambda,\mu}$ and $I_{\lambda,\mu}(|u|)=c_{\lambda,\mu}$. According to the proof of Theorem 4.3 of \cite{Willem1996book}, we can prove that  $I'_{\lambda,\mu}(|u|)=0$. Without loss of generality, we can assume $u\geq0$. By the theory of
elliptic regularity, $u\in C^2(\mathbb{R}^3)$, and by using strong maximum principle, we get $u>0$ in $\mathbb{R}^3$.
\hfill{$\Box$}
\vskip8pt

{\bf The proof of Theorem 1.1:}
By Lemma \ref{le:clambdaachieved}, there exists
$u_n\in \mathcal{N}_{\lambda_n,\mu}$  such that
$I_{\lambda_n,\mu}(u_n)=c_{\lambda_n,\mu}$,  $I'_{\lambda_n,\mu}(u_n)=0$ and $u_n(x)>0$, $x\in\mathbb{R}^3$. By  Lemma \ref{le:psbounded} and Lemma \ref{cbounded}, there exists $C>0$, such that $\|u_n\|_{\lambda_n}\leq C.$  Thus, noting \eqref{firstne}, there exists $u\in H^1(\mathbb{R}^3)$, such that
\begin{equation*}
\begin{aligned}
&u_n\rightharpoonup u\quad\hbox{in}\; H^1(\mathbb{R}^3),\\
&u_n\to u\quad\hbox{a.e. in}\  \mathbb{R}^3,\\
&u_n\rightharpoonup u
\quad\hbox{in}\  L^{q}(\mathbb{R}^3), \ \hbox{for }\ 2\leq q \leq 6.
\end{aligned}
\end{equation*}
\par
{\bf Step 1:} $u|_{\Omega^c}=0$, where $\Omega^c:=\{x\ |\ x\in\mathbb{R}^3\setminus \Omega\}.$
\par
If not, we have $u|_{\Omega^c}\neq 0.$
Then there exists a compact subset $F\subset \Omega^c$ with
$dist\{F, \partial\Omega\}>0$ such that $u|_F\neq 0$ and
\[\int_{F}u^2_ndx\to \int_F u^2dx>0.\]
Moreover, there exists
$\epsilon_0>0$ such that $V(x)\geq \epsilon_0$ for any $x\in F.$
\par
By the choice of $u_n$, we have
\begin{equation}\label{5lim1}
\displaystyle\|u_n\|^2_{\lambda_n}+F(u_n)=\mu\int_{\mathbb{R}^3}|u_n|^{p}dx+\int_{\mathbb{R}^3}|u_n|^{6}dx.
\end{equation}
Hence
\begin{equation*}
\begin{aligned}
I_{\lambda_n,\mu}(u_n)&=\displaystyle\frac{1}{2}\|u_n\|^2_{\lambda_n}+\frac{1}{4}F(u_n)-\frac{\mu}{p}\int_{\mathbb{R}^3}|u_n|^{p}dx-\frac{1}{6}\int_{\mathbb{R}^3}|u_n|^{6}dx\\
&\geq\displaystyle\frac{1}{2}\|u_n\|^2_{\lambda_n}+\frac{1}{4}F(u_n)-\frac{\mu}{p}\int_{\mathbb{R}^3}|u_n|^{p}dx-\frac{1}{p}\int_{\mathbb{R}^3}|u_n|^{6}dx\\
&=\displaystyle\frac{1}{2}\|u_n\|^2_{\lambda_n}+\frac{1}{4}F(u_n)-\frac{1}{p}(\|u_n\|^2_{\lambda_n}+F(u_n))\\
&\geq (\displaystyle\frac{1}{2}-\frac{1}{p})\int_{F}\lambda_nV(x)u^2_ndx\\
&\geq (\displaystyle\frac{1}{2}-\frac{1}{p})\lambda_n\epsilon_0\int_{F}u^2_ndx.
\end{aligned}
\end{equation*}
Then $I_{\lambda_n,\mu}(u_n)\to\infty$, as $n\to\infty$. This contradicts Lemma \ref{le:c-level}. Thus $u|_{\Omega^c}=0$ and so does
$u$ and $u \in H_0^1(\Omega)$.
\par
{\bf Step 2:} $u\neq0.$
\par
If not, by Lemma \ref{prop:strong}, we have $\int_{\mathbb{R}^3}|u_n|^{p}dx=o_n(1)$ and $F(u_n)=o_n(1)$. It follows from \eqref{5lim1} that
\[\|u_n\|^2_{\lambda_n}=\int_{\mathbb{R}^3}|u_n|^{6}dx+o_n(1).\]
Assume that
\[\|u_n\|^2_{\lambda_n}\to b\quad\text{and}\quad \int_{\mathbb{R}^3}|u_n|^{6}dx\to b.\]
By Lemma \ref{le:lowboundedinN}, $b\neq 0$. In term of the definition of $S$, we can get
\[S(\int_{\mathbb{R}^3}|u_n|^{6}dx)^{\frac{1}{3}}\leq\|u_n\|^2_{\lambda_n}.\]
Letting $n\to \infty$, we have
\begin{equation*}
S^{\frac{3}{2}}\leq b.
\end{equation*}
Thus, it follows from Remark \ref{r1} that
\begin{equation*}
\begin{aligned}
\displaystyle\frac{1}{3}S^{\frac{3}{2}}
&>\displaystyle\lim_{n\to\infty} c_{\lambda_n,\mu}\\
&= \displaystyle\lim_{n\to\infty}I_{\lambda_n,\mu}(u_n)\\
&=\displaystyle\lim_{n\to\infty}\Big[\frac{1}{2}\|u_n\|^2_{\lambda_n}+\frac{1}{4}F(u_n)-\frac{\mu}{p}\int_{\mathbb{R}^3}|u_n|^{p}dx-\frac{1}{6}\int_{\mathbb{R}^3}|u_n|^{6}dx\Big]\\
&=\displaystyle(\frac{1}{2}-\frac{1}{6})b\\
&\geq \displaystyle\frac{1}{3}S^{\frac{3}{2}},
\end{aligned}
\end{equation*}
which is a contrdiction.
\par
{\bf Step 3:} $I'_\mu(u)=0.$
\par
Consider a test function $\varphi\in \mathcal{C}_0^\infty(\Omega).$ Then
by the choice of $u_n$, we have $\langle I'_{\lambda_{n},\mu}(u_n), \varphi\rangle=0$. Thus
\begin{equation*}
\begin{aligned}
0&=\displaystyle\int_{\mathbb{R}^3}(\nabla u_n\nabla\varphi+\lambda_nV(x)u_n\varphi) dx+\int_{\mathbb{R}^3}\phi_{u_n}u_n\varphi dx\\
&\quad+\displaystyle\mu\int_{\mathbb{R}^3}|u_n|^{p-2}u_n\varphi dx+\int_{\mathbb{R}^3}|u_n|^{4}u_n\varphi dx.
\end{aligned}
\end{equation*}
Since $V(x)=0$ in $\Omega$, and let $n\to\infty$ in the above expression,
we have
\begin{equation*}
0=\displaystyle\int_{\Omega}\nabla u\nabla\varphi dx+\int_{\Omega}\phi_{u}u\varphi dx+\mu\int_{\Omega}|u|^{p-2}u\varphi dx+\int_{\Omega}|u|^{4}u\varphi dx,
\end{equation*}
which means $\langle I'_\mu(u), \varphi\rangle=0.$ By the theory of
elliptic regularity, $u\in C^2(\Omega)$, and by using strong maximum principle, we get $u>0$ in $\Omega$.
\par
{\bf Step 4:} $I_\mu(u)=c_\mu$.
\par
From the discussion in {\bf Step 3}, we know $u\in N_\mu$, so $I_\mu(u)\geq c_\mu$. Then
\begin{equation*}
\begin{aligned}
c_\mu&\leq I_\mu(u)\\
&=I_\mu(u)-\displaystyle\frac{1}{4}\langle I'_\mu(u),u\rangle\\
&=\displaystyle\frac{1}{4}\int_{\Omega}|\nabla u|^2dx+(\frac{1}{4}-\frac{1}{p})\mu\int_{\Omega}|u|^{p}dx+(\frac{1}{4}-\frac{1}{6})\int_{\Omega}|u|^{6}dx\\
&\leq\displaystyle\liminf_{n\to\infty}\Big[\frac{1}{4}\int_{\Omega}|\nabla u_n|^2dx+(\frac{1}{4}-\frac{1}{p})\mu\int_{\Omega}|u_n|^{p}dx+(\frac{1}{4}-\frac{1}{6})\int_{\Omega}|u_n|^{6}dx\Big]\\
&\leq\displaystyle\liminf_{n\to\infty}\Big[\frac{1}{4}\int_{\mathbb{R}^3}|\nabla u_n|^2dx+(\frac{1}{4}-\frac{1}{p})\mu\int_{\mathbb{R}^3}|u_n|^{p}dx+(\frac{1}{4}-\frac{1}{6})\int_{\mathbb{R}^3}|u_n|^{6}dx\Big]\\
&\leq\displaystyle\liminf_{n\to\infty}[I_{\lambda_n,\mu}(u_n)-\displaystyle\frac{1}{4}\langle I'_{\lambda_n,\mu}(u_n),u_n\rangle]\\
&=\displaystyle\liminf_{n\to\infty}c_{\lambda_n,\mu}\leq c_\mu,
\end{aligned}
\end{equation*}
Therefore, $I_\mu(u)=c_\mu$. Moreover, from the above proof, we know that
\begin{equation}\label{cc}
\begin{array}{ll}
\displaystyle \lim_{n\to \infty}I_{\lambda_n,\mu}(u_n)=\lim_{n\to \infty}c_{\lambda_n,\mu}=c_\mu=I_\mu(u).
\end{array}
\end{equation}
\par
{\bf Step 5:} $u_n\rightarrow u$ in $H^1(\mathbb{R}^3)$.
\par
In fact, from the proof in {\bf Step 4}, up to a subsequence, we have
\[\int_{\mathbb{R}^3}|\nabla u_n|^2dx\to\int_{\mathbb{R}^3}|\nabla u|^2dx,\quad \int_{\mathbb{R}^3}|u_n|^{p}dx\to\int_{\mathbb{R}^3}|u|^{p}dx,\]
\[\int_{\mathbb{R}^3}|u_n|^{6}dx\to\int_{\mathbb{R}^3}|u|^{6}dx,\quad
\int_{\mathbb{R}^3}\lambda_nV(x)u_n^2dx\to 0.\]
Therefore,
\[\|u_n-u\|_{\lambda_n}\to 0.\]
It follows that
\[\|u_n-u\|\to 0.\]
\hfill{$\Box$

\section{Remarks on the mountain pass levels}

In order to investigate the properties of  $c_{\mu}$, we need to analyse a functional on $H_0^1(\Omega)$, that is
\begin{equation*}
I_{0}(u)=\frac{1}{2}\|u\|^2_0+\frac{1}{4}F_0(u)-\frac{1}{6}\int_{\Omega}|u|^{6}dx.
\end{equation*}
Then we can define
\begin{equation*}
\mathcal{N}_{0}=\{u\in H_0^1(\Omega)\setminus \{0\}\ |\  \langle I'_{0}(u), u\rangle =0\},
\end{equation*}
and
\begin{equation*}
c_{0}=\inf_{u\in \mathcal{N}_{0}}I_{0}(u)
\end{equation*}

\begin{remark}\label{r2}
For $c_0$, $I_0$ and $\mathcal{N}_0$, there are similar results obtained from
lemma \ref{le:lowboundedinN} to Lemma \ref{le:moutain-pass geometry}.
\end{remark}
\par
We also need the  space $H_{0,rad}^1(B_r(0)):=\{u\in H_{0}^1(B_r(0))\ |\ u \ \text{is radial about the origin point}\}$, whose norm is given by
\begin{equation*}
\|u\|_{0,r}:=\Big(\int_{B_r(0)}|\nabla u|^2dx)\Big)^{1/2}.
\end{equation*}
For $u\in H_{0,rad}^1(B_r(0))$, we introduce another new functional
\begin{equation*}
I_{\mu,r}(u)=\frac{1}{2}\|u\|^2_{0,r}+\frac{1}{4}F_{0,r}(u)-\frac{\mu}{p}\int_{B_r(0)}|u|^{p}dx-\frac{1}{6}\int_{B_r(0)}|u|^{6}dx,
\end{equation*}
where $F_{0,r}(u)=\frac{1}{4\pi}\iint\limits_{B_r(0)\times B_r(0)}\frac{u^2(x)u^2(y)}{|x-y|}dxdy$. Similarly, we define
\begin{equation*}
\mathcal{N}_{\mu,r}=\{u\in H_{0,rad}^1(B_r(0))\setminus \{0\}\ |\  \langle I'_{\mu,r}(u), u\rangle =0\}
\end{equation*}
and
\begin{equation*}
c_{\mu,r}=\inf_{u\in \mathcal{N}_{\mu,r}}I_{\mu,r}(u).
\end{equation*}

\begin{remark}\label{r3}
For $c_{\mu,r}$, $I_{\mu,r}$ and $\mathcal{N}_{\mu,r}$, there are similar results obtained from
lemma \ref{le:lowboundedinN} to Lemma \ref{le:c-level}. By Mountain pass Theorem, we can see that there exists $u\in H_{0,rad}^1(B_r(0))$ such that
$I_{\mu,r}(u)=c_{\mu}$ and $I'_{\mu,r}(u)=0$.
\end{remark}

\begin{lemma}\label{le:00limit}
$\displaystyle\lim_{\mu\to 0}c_{\mu}=c_0.$
\end{lemma}
{\bf Proof:} We just need to prove that, for any $\mu_n\to 0$, as $n\to \infty$, we have  \[\lim_{n\to \infty}c_{\mu_n}=c_0.\]
In fact, by the definition of $c_{\mu_n}$ and $c_0$, we can obtain $c_{\mu_n}\leq c_0$. Then
\begin{equation}\label{lim1}
\limsup_{n\to \infty}c_{\mu_n}\leq c_0.
\end{equation}
On the other hand, by Remark \ref{ro}, we have $u_n\in H_0^1(\Omega)$ such that
\[I_{\mu_n}(u_n)=c_{\mu_n} \quad\text{and}\quad I'_{\mu_n}(u_n)=0.\]
It follows from Remark \ref{ro} that there exists $t_n>0$ such that $t_nu_n\in N_0$. Then
\[c_0\leq I_{0}(t_nu_n)=I_{\mu_n}(t_nu_n)+\frac{\mu_nt_n^p}{p}\int_{\mathbb{R}^3}|u_n|^pdx,\]
hence
\[c_0\leq c_{\mu_n}+\frac{\mu_nt_n^p}{p}|u_n|^p.\]
If $t_n$ is bounded, then we have
\begin{equation}\label{lim2}
c_0\leq \liminf_{n\to \infty}c_{\mu_n}.
\end{equation}
The result follows \eqref{lim1} and \eqref{lim2}.

Now we prove the boundness of $\{t_n\}$. Assume by contradiction that $t_n\to \infty$. By the choice of $u_n$ and $c_{\mu_n}<\frac{1}{3}S^{\frac{3}{2}}$, we know $\|u_n\|_0$ is bounded. Since $t_nu_n\in N_0$, we have
\begin{equation*}
\frac{1}{t_n^4}\|u_n\|^2_0+\frac{1}{t_n^2}F_0(u_n)=\int_{\Omega}|u_n|^{6}dx.
\end{equation*}
Thus
\begin{equation}\label{lim3}
\int_{\Omega}|u_n|^{6}dx\to 0,\quad \text{as}\ n\to\infty.
\end{equation}
By $I'_{\mu_n}(u_n)=0$, we can get
\begin{equation*}
\|u_n\|^2_0+F_0(u_n)=\mu_n\int_{\Omega}|u_n|^{p}dx+\int_{\Omega}|u_n|^{6}dx.
\end{equation*}
Then it follows from \eqref{lim3} that  $\|u_n\|^2_0=o_n(1)$. Therefore, $c_{\mu_n}\to 0$, which is absurd.
\hfill{$\Box$}

%引理4.4

\begin{lemma}\label{le:0limit}
$\displaystyle c_0=\frac{1}{3}S^{3/2}.$
\end{lemma}
{\bf Proof:} Similar to the proof of Lemma \ref{le:c-level}, we can get
\[c_0\leq \frac{1}{3}S^{3/2}+o_\epsilon(1).\]
Passing to the limit yields $c_0\leq \frac{1}{3}S^{3/2}$.

On the other hand, there exists $\{u_n\}\subset H_0^1(\Omega)$ such that
\[I_{0}(u_n)\to c_{0} \quad\text{and}\quad I'_{0}(u_n)\to 0.\]
By $I'_{0}(u_n)\to 0$,  we have
\[o_n(1)=\langle I'_{0}(u_n), u_n\rangle=\|u_n\|^2_0+F_0(u_n)-\int_{\Omega}|u_n|^{6}dx.\]
Assume that
\[\|u_n\|^2_0\to l_1,\quad F_0(u_n)\to l_2 \ \text{and}\ \int_{\Omega}|u_n|^{6}dx\to l_3.\]
Then,  we have
\begin{equation}\label{1equality2}
l_1+l_2=l_3.
\end{equation}
By the definition of $S$, we can get
\[S(\int_{\Omega}|u_n|^{6}dx)^{\frac{1}{3}}\leq\|u_n\|^2_0.\]
Letting $n\to \infty$, we have
\begin{equation}\label{2equality2}
Sl_2^{\frac{1}{3}}\leq l_1.
\end{equation}
From \eqref{1equality2} and \eqref{2equality2}, wo know that
\begin{equation}\label{3equality2}
l_1\geq S^{3/2}\quad \text{and} \quad l_3\geq S^{3/2}.
\end{equation}
By $I_{0}(u_n)\to c_{0}$, we have
\begin{equation*}
\begin{aligned}
c_0&=I_0(u_n)+o_n(1)\\
&=\displaystyle\frac{1}{2}\|u_n\|^2_0+\frac{1}{4}F_0(u_n)-\frac{1}{6}\int_{\Omega}|u_n|^{6}dx+o_n(1)\\
&=\displaystyle\frac{1}{2}\|u_n\|^2_0-\frac{1}{4}(\|u_n\|^2_0-\int_{\Omega}|u_n|^{6}dx)-\frac{1}{6}\int_{\Omega}|u_n|^{6}dx+o_n(1)\\
&=\displaystyle\frac{1}{4}\|u_n\|^2_0+\frac{1}{12}\int_{\Omega}|u_n|^{6}dx+o_n(1)\\
&=\displaystyle\frac{1}{4}l_1+\frac{1}{12}l_3+o_n(1).
\end{aligned}
\end{equation*}
Noting \eqref{3equality2}, it follows from the above expression that
$c_0\geq \frac{1}{3}S^{3/2}.$
\hfill{$\Box$}

%引理4.5

\begin{lemma}\label{le:limit}
For any $r>0$, $\displaystyle\lim_{\mu\to 0}c_{\mu}=\lim_{\mu\to 0}c_{\mu,r}=\frac{1}{3}S^{3/2}.$
\end{lemma}
{\bf Proof:}  It follows from Lemma \ref{le:00limit} and Lemma \ref{le:0limit} that
\[\lim_{\mu\to 0}c_{\mu}=\frac{1}{3}S^{3/2}.\]
According to the definition of $c_{\mu,r}$, we have
\begin{equation*}
c_{\mu,r}=\inf_{u\in N_{\mu,r}}I_{\mu,r}(u)\geq \inf_{\{u\in H_0^1(B_r(0))\ |\ \langle I'_{\mu,r}(u), u\rangle=0\}}I_{\mu,r}(u).
\end{equation*}
Since the limit of the last term in the above expression is $c_\mu$,  we can get
\begin{equation*}
\liminf_{\mu\to 0}c_{\mu,r}\geq \frac{1}{3}S^{3/2}.
\end{equation*}
On the other hand, similar to the proof of Lemma \ref{le:c-level}, and choosing $\varphi$ to be a radial function, we can get
\begin{equation*}
c_{\mu,r}< \frac{1}{3}S^{3/2}.
\end{equation*}
Then we can obtain
\begin{equation*}
\limsup_{\mu\to 0}c_{\mu,r}\leq \frac{1}{3}S^{3/2}.
\end{equation*}
Therefore, for any $r>0$, $\displaystyle\lim_{\mu\to 0}c_{\mu,r}=\frac{1}{3}S^{3/2}.$
\hfill{$\Box$}

Since $\Omega$ is a smooth bounded domain, we can fix $r>0$ small enough such that
\begin{equation*}
\Omega^+_r=\{ x\in\mathbb{R}^3 \ |\  dist(x,\Omega)\leq r \}
\end{equation*}
and
\begin{equation*}
\Omega^-_r=\{ x\in\Omega \ |\  dist(x,\partial\Omega)\geq r \}
\end{equation*}
are homotopically equivalent to $\Omega$. Moreover, we may assume that $B_r(0)\subset \Omega$.

For $0\neq u\in L^6(\Omega)$, we consider its center of mass
\begin{equation*}
\beta_0(u):=\frac{\int_{\Omega} xu^6dx}{\int_{\Omega} u^6dx }
\end{equation*}

\begin{lemma}\label{le:importance} There exists $\mu^*>0$ such that if $\mu\in (0,\mu^*)$ and  $u\in \mathcal{N}_\mu$ with $I_\mu(u)\leq c_{\mu,r}+o_\mu(1)$, then $\beta_0(u)\in \Omega^+_{r/2}$.
\end{lemma}
{\bf Proof:} Suppose by contradiction that there exist $\mu_n\to 0$, $u_n\in \mathcal{N}_{\mu_n}$ and
$I_{\mu_n}(u_n)\leq c_{\mu_n,r}+o_n(1)$ such that $\beta_0(u_n)\notin \Omega^+_{r/2}$. From $u_n\in \mathcal{N}_{\mu_n}$ and $I_{\mu_n}(u_n)\leq c_{\mu_n,r}+o_n(1)$, it is easy to obtain that $\{u_n\}$ is bounded in $H_0^1(\Omega)$. It follows from  $u_n\in \mathcal{N}_{\mu_n}$ that
\begin{equation}\label{equality0}
\|u_n\|^2_0+F_0(u_n)=\mu_n\int_{\Omega}|u_n|^{p}dx+\int_{\Omega}|u_n|^{6}dx.
\end{equation}
Thus, noting $\mu_n\to 0$,
\begin{equation*}
\begin{aligned}
c_{\mu_n}\leq I_{\mu_n}(u_n)&=\displaystyle\frac{1}{2}\|u_n\|^2_0+\frac{1}{4}F_0(u_n)-\frac{\mu_n}{p}\int_{\Omega}|u_n|^{p}dx-\frac{1}{6}\int_{\Omega}|u_n|^{6}dx\\
&=\displaystyle\frac{1}{4}\|u_n\|^2_0+(\frac{1}{4}-\frac{1}{p})\mu_n\int_{\Omega}|u_n|^{p}dx+\frac{1}{12}\int_{\Omega}|u_n|^{6}dx\\
&=\displaystyle\frac{1}{4}\|u_n\|^2_0+\frac{1}{12}\int_{\Omega}|u_n|^{6}dx+o_n(1)\\
&\leq c_{\mu_n,r}+o_n(1)
\end{aligned}
\end{equation*}
Assume that $\|u_n\|^2_0\to l_1$ and $\int_{\Omega}|u_n|^{6}dx\to l_2$. Then, by Lemma \ref{le:limit}, we have
\begin{equation}\label{equality1}
\frac{1}{4}l_1+\frac{1}{12}l_2=\frac{1}{3}S^{3/2}.
\end{equation}
By the definition of $S$, we can get
\[S(\int_{\Omega}|u_n|^{6}dx)^{\frac{1}{3}}\leq\|u_n\|^2_0.\]
Letting $n\to \infty$, we have
\begin{equation}\label{equality2}
Sl_2^{\frac{1}{3}}\leq l_1.
\end{equation}
From \eqref{equality0}, wo know that
\begin{equation*}
\|u_n\|^2_0\leq\int_{\Omega}|u_n|^{6}dx+o_n(1).
\end{equation*}
Then,
\begin{equation}\label{equality3}
l_1\leq l_2.
\end{equation}
It follows from \eqref{equality1}, \eqref{equality2} and \eqref{equality3} that $l_2= l_1= S^{3/2}$.
\par
Define
\[\omega_n=\displaystyle\frac{u_n}{\int_{\Omega}|u_n|^{6}dx}.\]
Then we have
\[\int_{\Omega}|\omega_n|^{6}dx=1\]
and
\begin{equation*}
\int_{\Omega}|\nabla\omega_n|^{6}dx\to S,\quad \text{as}\  n\to\infty.
\end{equation*}
Using the method in \cite{MR1974041}, we can assume $\omega_n$ is nonnegative. By Lemma 3.1 in \cite{MR1974041}, there exists  $(y_n,\lambda_n)\in\mathbb{R}^N\times\mathbb{R}^+$ with $\lambda_n\to 0$ and $y_n\to y\in \bar{\Omega}$ such that
\begin{equation*}
v_n(x):=\lambda_n^{\frac{1}{2}}\omega_n(\lambda_nx+y_n)\to\omega \ \ \text{in}\ D^{1,2}(\mathbb{R}^3),\quad \text{as}\  n\to\infty.
\end{equation*}
Suppose that $\phi\in C_0^\infty(\mathbb{R}^3)$ satisties $\phi(x)=x$, $x\in \Omega$. Then
\begin{equation*}
\beta_0(u_n)=\int_{\Omega}\phi(x)\omega_n^6(x)dx=\int\phi(\lambda_nx+y_n)v_n^6(x)dx
\end{equation*}
Now, Lebesgue Dominated Convergence Theorem yields
\begin{equation*}
\beta_0(u_n)\to y\in \bar{\Omega}, \quad \text{as} \ n\to \infty,
\end{equation*}
which contradicts our assumption that $\beta_0(u_n)\notin \Omega^+_{r/2}$.
\hfill{$\Box$}

As in \cite{MR1762697}, we choose $R>2\text{diam}(\Omega)$ with $\Omega\subset B_R(0)$ and set
\begin{equation*}
\psi(t)=
\begin{cases}
1, \quad 0\leq t\leq R,\\
\frac{R}{t}, \quad t\geq R.\\
\end{cases}
\end{equation*}
Define
\begin{equation*}
\beta(u)=\frac{\int_{\mathbb{R}^3}\psi(|x|)|u|^6xdx}{\int_{\mathbb{R}^3}|u|^6dx},\quad \text{for}\ u\in L^6(\mathbb{R}^3)\setminus\{0\}.
\end{equation*}

\begin{lemma}\label{6importance}
There is $\lambda^*>0$ and $\mu^*>0$  such that if $u\in \mathcal{N}_{\lambda,\mu}$ and
$I_{\lambda,\mu}(u)\leq c_{\mu,r}$, then  $\beta(u)\in\Omega^+_r$ for all $\lambda\geq\lambda^*$ and $\mu<\mu^*$.
\end{lemma}
{\bf Proof:} Assume, by contradiction, that for $\mu$ arbitrarily small there exist  $\lambda_n\to\infty$ and $u_n\in \mathcal{N}_{\lambda_n,\mu}$ with $I_{\lambda_n,\mu}(u_n)\leq c_{\mu,r}$ and
\[\beta(u_n)\notin \Omega^+_r.\]
It is easy to see  that $\{\|u_n\|_{\lambda_n}\}$ is bounded. From {\bf Step 1} in the proof of Theorem $1.1$, there exists $u\in H^1_0(\Omega)$ such that
\begin{equation*}
\begin{aligned}
&u_n\rightharpoonup u\quad\hbox{in}\; H^1(\mathbb{R}^N),\\
&u_n\to u\quad\hbox{a.e. in}\  \mathbb{R}^3,\\
&u_n\rightharpoonup u
\quad\hbox{in}\  L^{q}(\mathbb{R}^3), \ \hbox{for }\ 2\leq q \leq 6.
\end{aligned}
\end{equation*}
Moreover, it follows from Lemma \ref{prop:strong} that
\begin{equation*}
u_n\to u\quad\hbox{in}\  L^{q}(\mathbb{R}^3), \ \hbox{for }\ 2< q < 6.
\end{equation*}
Let $v_n=u_n-u$. Then, noting $u_n\in \mathcal{N}_{\lambda_n,\mu}$, we have
\begin{equation}\label{6lim1}
\begin{aligned}
\|v_n\|^2_{\lambda_n}&=\displaystyle\|u_n\|^2_{\lambda_n}-\int_{\mathbb{R}^3}|\nabla u|^{2}dx+o_n(1)\\
&=\displaystyle\mu\int_{\mathbb{R}^3}|u_n|^{p}dx+\int_{\mathbb{R}^3}|u_n|^{6}dx-F(u_n)
-\int_{\mathbb{R}^3}|\nabla u|^{2}dx+o_n(1)\\
&=\displaystyle\mu\int_{\mathbb{R}^3}|v_n|^{p}dx+\int_{\mathbb{R}^3}|v_n|^{6}dx-F(v_n)\\
&\quad\displaystyle+\mu\int_{\mathbb{R}^3}|u|^{p}dx+\int_{\mathbb{R}^3}|u|^{6}dx-F(u)-\int_{\mathbb{R}^3}|\nabla u|^{2}dx+o_n(1)\\
&=\displaystyle\int_{\mathbb{R}^3}|v_n|^{6}dx+\mu\int_{\mathbb{R}^3}|u|^{p}dx+\int_{\mathbb{R}^3}|u|^{6}dx-F(u)-\int_{\mathbb{R}^3}|\nabla u|^{2}dx+o_n(1).
\end{aligned}
\end{equation}
\par
Next, we divide into two cases.\\
{\bf Case 1:}  $\displaystyle\mu\int_{\mathbb{R}^3}|u|^{p}dx+\int_{\mathbb{R}^3}|u|^{6}dx\leq F(u)+\int_{\mathbb{R}^3}|\nabla u|^{2}dx$.
\par
From \eqref{6lim1}, we have
\begin{equation}\label{6lim2}
\|v_n\|^2_{\lambda_n}\leq \int_{\mathbb{R}^3}|v_n|^{6}dx+o_n(1).
\end{equation}
We claim that $\displaystyle\int_{\mathbb{R}^3}|v_n|^{6}dx\to 0$. If not, there exists $b>0$ such that
\begin{equation}\label{6lim3}
\int_{\mathbb{R}^3}|v_n|^{6}dx\to b.
\end{equation}
By the definition of $S$ and \eqref{6lim2}, we have
\begin{equation}
S(\int_{\mathbb{R}^3}|v_n|^{6}dx)^{\frac{1}{3}}\leq\|v_n\|^2_{\lambda_n}\leq\int_{\mathbb{R}^3}|v_n|^{6}dx+o_n(1)
\end{equation}
Thus, it follows from \eqref{6lim3} that $b\geq S^{\frac{3}{2}}$.
\par
On the other hand, noting $u_n\in \mathcal{N}_{\lambda_n,\mu}$, we have
\begin{equation*}
\begin{aligned}
I_{\lambda_n,\mu}(u_n)&=\displaystyle\frac{1}{2}\|u_n\|^2_{\lambda_n}+\frac{1}{4}F(u_n)-\frac{\mu}{p}\int_{\mathbb{R}^3}|u_n|^{p}dx-\frac{1}{6}\int_{\mathbb{R}^3}|u_n|^{6}dx\\
&=\displaystyle\frac{1}{4}\|u_n\|^2_{\lambda_n}+(\frac{\mu}{4}-\frac{\mu}{p})\int_{\mathbb{R}^3}|u_n|^{p}dx+\frac{1}{12}\int_{\mathbb{R}^3}|u_n|^{6}dx\\
&=\displaystyle\frac{1}{4}\|v_n\|^2_{\lambda_n}+\int_{\mathbb{R}^3}|\nabla u|^{2}dx+(\frac{\mu}{4}-\frac{\mu}{p})\int_{\mathbb{R}^3}|u|^{p}dx+\frac{1}{12}\int_{\mathbb{R}^3}|v_n|^{6}dx\\
&\quad\displaystyle+\frac{1}{12}\int_{\mathbb{R}^3}|u|^{6}dx+o_n(1)\\
&\geq\displaystyle\frac{1}{4}S(\int_{\mathbb{R}^3}|v_n|^{6}dx)^{\frac{1}{3}}+\frac{1}{12}\int_{\mathbb{R}^3}|u|^{6}dx+o_n(1).
\end{aligned}
\end{equation*}
Then, it follows from $I_{\lambda_n,\mu}(u_n)\leq c_{\mu,r}$ and \eqref{6lim3} that
\begin{equation*}
\frac{1}{4}Sb^{\frac{1}{3}}+\frac{1}{12}b\leq c_{\mu,r}.
\end{equation*}
Thus, we have $c_{\mu,r}\geq\frac{1}{3}S^{\frac{3}{2}}$, which is absurd. Consequently, $u_n\to u$ in $L^6(\mathbb{R}^3)$ and, therefore, $\beta(u_n)\to\beta(u)=\beta_0(u)$. By \eqref{6lim2}, we have $u\in \mathcal{N}_\mu$ and $I_{\mu}(u)=\lim_{n\to\infty}I_{\lambda_n,\mu}(u_n)\leq c_{\mu,r}$. Then it follows from Lemma \ref{le:importance} that $\beta(u)\in \Omega^+_{r/2}$. This contradicts our assumption that $\beta(u_n)\notin \Omega^+_r$.\\
{\bf Case 2:} $\displaystyle\mu\int_{\mathbb{R}^3}|u|^{p}dx+\int_{\mathbb{R}^3}|u|^{6}dx> F(u)+\int_{\mathbb{R}^3}|\nabla u|^{2}dx$.
\par
In this case, there exists $t_\mu\in(0,1)$ such that $t_\mu u\in \mathcal{N}_\mu$ and, therefore,
\begin{equation*}
\begin{aligned}
c_\mu &\leq I_{\mu}(t_\mu u)\\
&=\displaystyle\frac{t_\mu^2}{2}\|u\|^2_0+\frac{t_\mu^4}{4}F(u)-\frac{t_\mu^p}{p}\mu\int_{\mathbb{R}^3}|u|^{p}dx-\frac{t_\mu^6}{6}\int_{\mathbb{R}^3}|u|^{6}dx\\
&=\displaystyle\frac{t_\mu^2}{3}\|u\|^2_0+\frac{t_\mu^4}{12}F(u)-(\frac{1}{p}-\frac{1}{6})t_\mu^p\mu\int_{\mathbb{R}^3}|u|^{p}dx.
\end{aligned}
\end{equation*}
It follows from  $u_n\in \mathcal{N}_{\lambda_n,\mu}$ that
\begin{equation*}
\begin{aligned}
I_{\lambda_n,\mu}(u_n)&=\displaystyle\frac{1}{2}\|u_n\|^2_{\lambda_n}+\frac{1}{4}F(u_n)-\frac{\mu}{p}\int_{\mathbb{R}^3}|u_n|^{p}dx-\frac{1}{6}\int_{\mathbb{R}^3}|u_n|^{6}dx\\
&=\displaystyle\frac{1}{3}\|u_n\|^2_{\lambda_n}+\frac{1}{12}F(u_n)-(\frac{1}{p}-\frac{1}{6})\mu\int_{\mathbb{R}^3}|u_n|^{p}dx.
\end{aligned}
\end{equation*}
Thus,
\begin{equation*}
\begin{aligned}
\displaystyle c_\mu+(\frac{1}{p}-\frac{1}{6})t_\mu^p\mu\int_{\mathbb{R}^3}|u|^{p}dx&\leq
\displaystyle\frac{t_\mu^2}{3}\|u\|^2_0+\frac{t_\mu^4}{12}F(u)\\
&\leq\displaystyle\liminf_{n\to\infty}\Big(\frac{t_\mu^2}{3}\|u_n\|^2_0+\frac{t_\mu^4}{12}F(u_n)\Big)\\
&\leq\displaystyle\liminf_{n\to\infty}\Big(\frac{1}{3}\|u_n\|^2_0+\frac{1}{12}F(u_n)\Big)\\
&\leq\displaystyle\liminf_{n\to\infty}\Big(I_{\lambda_n,\mu}(u_n)+(\frac{1}{p}-\frac{1}{6})\mu\int_{\mathbb{R}^3}|u_n|^{p}dx\Big)\\
&\leq \displaystyle c_{\mu,r}+(\frac{1}{p}-\frac{1}{6})\mu\int_{\mathbb{R}^3}|u|^{p}dx.
\end{aligned}
\end{equation*}
Then, by Lemma \ref{le:limit}, we have $t_\mu=1+o_\mu(1)$. Consequently,
\begin{equation}\label{6lim5}
\Big|\int_{\mathbb{R}^3}|\nabla u_n|^2dx-\int_{\mathbb{R}^3}|\nabla (t_\mu u)|^2dx\Big|\leq 3(c_{\mu,r}-c_\mu)+o_\mu(1)+o_n(1).
\end{equation}
Noting  $u_n\rightharpoonup u$ in $H^1(\mathbb{R}^3)$, it follows from \eqref{6lim5}
that
\begin{equation*}
\int_{\mathbb{R}^3}|\nabla (u_n-t_\mu u)|^2dx\leq 3(c_{\mu,r}-c_\mu)+o_\mu(1)+o_n(1).
\end{equation*}
Therefore,
\begin{equation*}
\int_{\mathbb{R}^3}(u_n-t_\mu u)^6dx\leq \frac{1}{S^3}\Big[3(c_{\mu,r}-c_\mu)+o_\mu(1)+o_n(1)\Big]^3.
\end{equation*}
Then, for $\mu$ small and $n$ large, we have
\begin{equation*}
|\beta(u_n)-\beta(t_\mu u)|\leq \frac{r}{2}.
\end{equation*}
From $I_{\lambda_n,\mu}(u_n)\leq c_{\mu,r}$, it is easy to see that $I_{\mu}(t_\mu u)\leq c_{\mu,r}+o_\mu(1)$. Thus, by Lemma \ref{le:importance}, we can obtain $\beta(t_\mu u)=\beta_0(t_\mu u)\in \Omega^+_{r/2}$. This contradicts our assumption that $\beta(u_n)\notin \Omega^+_r$.
\hfill{$\Box$

\section{The existence of multiple solutions}

\begin{lemma}\label{le5.1}
There are $\lambda^*>0$ large and $\mu^*>0$ small such that
\[c_{\mu,r}<2c_{\lambda,\mu},\quad \text{for any}\  \lambda\geq\lambda^* \text{and}\  \mu\in(0,\mu^*)\]
\end{lemma}
{\bf Proof:} It is easy to get the result from \eqref{cc} and Lemma \ref{le:limit}.
\hfill{$\Box$

\begin{lemma}\label{le5.2}
Assume that $c_{\mu,r}<2c_{\lambda,\mu}$ and $u\in E$ is a nontrivial critical point of $I_{\lambda,\mu}$ with $I_{\lambda,\mu}(u)\leq c_{\mu,r}$. Then, $u$ is positive or  $u$ is negative.
\end{lemma}
{\bf Proof:} If $u^\pm\neq 0$, then $u^\pm\in \mathcal{N}_{\lambda,\mu}$, and so,
\[c_{\mu,r}\geq I_{\lambda,\mu}(u)=I_{\lambda,\mu}(u^+)+I_{\lambda,\mu}(u^-)\geq 2c_{\lambda,\mu}\]
which is absurd.
\hfill{$\Box$

\begin{remark}\label{positiveso}
Since $I_{\lambda,\mu}$ is even, we can assume that nontrivial critical point is positive.
\end{remark}

\begin{lemma}\cite{MR1196690}\label{chang}  Let $M$ be a $C^1$ functional defined on a $C^1$ Finsler manifold $M.$ If $I$ is bounded from below and satisfies the $(PS)$ condition, then $I$ has at least $cat_MM$ distinct critical points.
\end{lemma}

\begin{lemma}\cite{MR1088278}\label{cerami2}  Let $\Gamma,\Omega^+,\Omega^-$ be closed sets with $\Omega^-\subset\Omega^+.$ Let $\Phi:\Omega^-\rightarrow\Gamma,$ $\beta:\Gamma\rightarrow\Omega^+$ be two continuous maps such that $\beta\circ\Phi$ is homotopically equivalent to the embedding $Id:\Omega^-\rightarrow\Omega^+.$ Then $cat_\Gamma\Gamma\geq cat_{\Omega^+}\Omega^-.$
\end{lemma}

In what follows, $u_r\in H^1_0(B_r(0))$ is a positive radial ground state solution for the functional $I_{\mu,r}$, that is
\begin{equation*}
 I_{\mu,r}(u_r)=c_{\mu,r}\quad\text{and}\quad I'_{\mu,r}(u_r)=0.
\end{equation*}
Define operator $\Psi_r\ :\ \Omega^-_r\to H^1_0(\Omega)$ by
\begin{equation*}
\Psi_r(y)(x)=
\begin{cases}
u_r(|x-y|),\quad x\in B_r(y),\\
0, \quad x\in \Omega\setminus B_r(y)
\end{cases}
\end{equation*}
which is continuous and satisfies
\begin{equation}\label{b1thm}
\beta(\Psi_r(y))=y,\quad y\in \Omega^-_r.
\end{equation}
Moreover,
\begin{equation}\label{b2thm}
\Psi_r(y)(x)\in \mathcal{N}_{\lambda,\mu}\quad\text{and}\quad I_{\lambda,\mu}(\Psi_r(y)(x))=I_{\mu,r}(\Psi_r(y)(x))=c_{\mu,r}.
\end{equation}
\\
{\bf The proof of Theorem 1.2:}
\par
For $\lambda>\lambda^*$ and $\mu<\mu^*$, we define two maps
\begin{equation*}
\Omega^-_r \xrightarrow{\Psi_r} I^{c_{\mu,r}}_{\lambda,\mu} \xrightarrow{\beta}\Omega^+_r.
\end{equation*}
where $I^{c_{\mu,r}}_{\lambda,\mu}=\{u\in \mathcal{N}_{\lambda,\mu}\ |\ I_{\lambda,\mu}\leq c_{\mu,r}\}$.
By Lemma \ref{6importance} and \eqref{b2thm}, the maps are well defined.
Since $I_{\lambda,\mu}$ satisfies the $(PS)_c$ condition on $\mathcal{N}_{\lambda,\mu}$ for $c\leq c_{\mu,r}$. By Lemma \ref{chang}, $I_{\lambda,\mu}$ has at least $cat_{I^{c_{\mu,r}}_{\lambda,\mu}} (I^{c_{\mu,r}}_{\lambda,\mu})$
critical points on $\mathcal{N}_{\lambda,\mu}$. Then, \eqref{b1thm} and Lemma \ref{cerami2} ensures that $I_{\lambda,\mu}$ has at least  $cat_{\Omega^+_r} (\Omega^-_r)$ critical points on $\mathcal{N}_{\lambda,\mu}$, and consequently, $cat_{\overline{\Omega}} (\overline{\Omega})$ critical points in $E$. Therefore, by Lemma \ref{le5.1}, Lemma \ref{le5.2} and Remark \ref{positiveso}, problem \eqref{maineq} has at least $cat_{\overline{\Omega}} (\overline{\Omega})$ positive solutions.

%%%%%%%%%%%%%%%%%%%%%%%%%%%%%%

\section*{Acknowledgments}
We would like to thank the anonymous referee for his/her careful readings of our manuscript and the useful comments made for its improvement. The first author thanks his advisor Prof. Zhongwei Tang for suggestions and help. The second author also thanks the support of RTG 2419 by the German Science Foundation (DFG).

\bibliographystyle{plain}
\bibliography{system1}
\end{document}